\documentclass[a4paper,11pt]{article}
\usepackage{latexsym}
\usepackage{amsmath}
\usepackage{amssymb}
\usepackage{bm}
\usepackage{graphicx}
\usepackage{wrapfig}
\usepackage{fancybox}
\begin{document}
\begin{center}
{\bf GROUP EMBEDDING INTO WREATH PRODUCTS: CARTESIAN DECOMPOSITION APPROACH}

{\bf Enoch Suleiman}

Department of Mathematics, Federal University Gashua, Yobe State, Nigeria.

enochsuleiman@gmail.com

{\bf Muhammed Salihu Audu}

Department of Mathematics, University of Jos, Plateau State, Nigeria.

audumso2@gmail.com

{\bf Sunday U. Momoh}

Department of Mathematics, University of Jos, Plateau State, Nigeria.

smomohu@yahoo.com
\end{center}
{\bf Abstract:} In this paper, we showed how a group acting regularly and a diagonal group are embedded into the wreath products in there product action using the Cartesian Decomposition.

{\bf Keywords:} Permutation Group, Embedding, Wreath Products, Cartesian Decomposition, product action.

{\bf 2020 Subject Classification:} 20B05,20E22, 20B35, 20B30, 20B99

{\bf 1. INTRODUCTION}

	Praeger and Schneider [1] proved that for a given $X$, a subgroup that is transitive of a given wreath product $ Sym\Gamma wrSym\Delta$ on $\Delta$, then $X$ is shown to be isomorphic to a subgroup of the wreath product of the permutation group prompted by its stabilizer $X_{\delta}$ on the set $\Gamma$ and also a given group prompted from a set $X$ on $\Delta$. Preager and Schneider [2] also proved that quasiprimitive permutation groups that is of simple diagonal type is not in any way isomorphic to a subgroup of wreath products that is acting on the same point set. 

	Many other people have made some progress on embedment of Groups into wreath products [3, 5, 6, 7,10]. In this paper, we proved that a group that is acting regularly on a given set and a diagonal group acting on a set in product action are embeddable into wreath products in such actions.

{\bf 2. NOTION OF CARTESIAN DECOMPOSITION AND WREATH PRODUCTS}

{\bf Definition 2.1} [9]{\bf :} Cartesian decomposition is define given a set $\Omega$ that is finite of partitions of $\Omega,\ \epsilon=\{\Gamma_{1},\Gamma_{2},\cdots,\Gamma_{k}\}$, with $|\Gamma_{i}|\geq 2$, for all $i\geq 1$ and $|\gamma_{1}\cap\gamma_{2}\cap\cdots\cap\gamma_{k}|=1$ for all $\gamma_{1}\in\Gamma_{1},\gamma_{2}\in\Gamma_{2},\cdots,\gamma_{k}\in\Gamma_{k}$.

	Any Cartesian decomposition is called a {\it trivial Cartesian decomposition} if it comprises of just a single partition, that is partition into singletons. Cartesian decomposition is called {\it homogeneous} if it has the property that all the $\Gamma_{i}$ have the equal number of elements. If $\{\Gamma_{1},\Gamma_{2},\cdots,\Gamma_{k}\}$, is a particular cartesian decomposition of a given set $\Omega$, then the defining property yields a well-defined one-to-one correspondence between $\Omega$ and $\Gamma_{1}\times\Gamma_{2}\times\cdots\times\Gamma_{k}$  given by $\omega\mapsto\left(\gamma_{1},\gamma_{2},\cdots,\gamma_{k}\right)$ where, for a given $i =1,2,\cdots,k$, the block $\gamma_{i}\in\Gamma_{i}$ is the unique block of $\Gamma_{i}$ which contains $\omega$. Thus the set $\Omega$ can be obviously recognized with the cartesian product $\Gamma_{1}\times\Gamma_{2}\times\cdots\times\Gamma_{k}$.

{\bf Example 2.2:} Let $\Omega=(1,2)^{3}$, we have the subsequent partitions of $\Omega$ as.
$$
\Gamma_{1}=\{\{(1,1,1),(1,1,2),(1,2,1),(1,2,2)\},\{(2,1,1),(2,1,2),(2,2,1),(2,2,2)\}\}
$$
$$
\Gamma_{2}=\{\{(1,1,1),(1,1,2),(2,1,1),(2,1,2)\},\{(1,2,1),(1,2,2),(2,2,1),(2,2,2)\}\}
$$
$$
\Gamma_{3}=\{\{(1,1,1),(1,2,1),(2,1,1),(2,2,1)\},\{(1,1,2),(1,2,2),(2,1,2),(2,2,2)\}\}
$$
With $|\gamma_{1}\cap\gamma_{2}\cap\gamma_{3}|=1$ for all $\gamma_{1}\in\Gamma_{1},\gamma_{1}\in\Gamma_{2},\gamma_{3}\in\Gamma_{3}$ and $|\gamma_{1}|=|\gamma_{2}|=|\gamma_{3}|=2$.

{\bf Definition 2.3(Wreath products of groups} [1,9] {\bf )}

Suppose $G$ and $H$ are two given groups with $\phi$ a homomorphism from $H$ to the automorphism group $Aut\left(G\right)$. The semidirect product of $G$ and $H$, denoted by $G\rtimes_{\phi}H$ or simply by $G\rtimes H$, is given as follows. The underlying set of the group $G\rtimes H$ is the direct product $G\times H$ of sets and the multiplication of two elements $\left(g_{1},h_{1}\right)$ and $\left(g_{2},h_{2}\right)$ is defined as 
$$
\left(g_{1},h_{1}\right)\left(g_{2},h_{2}\right):=\left(g_{1}\left(g_{2}\left(h_{1}^{-1}\phi\right)\right),h_{1}h_{2}\right)
$$
It is routine to check that the semidirect product $G\rtimes H$ is a group. I can easily be showned that $\overline{G}=\left\{\left(g,1\right)|g\in G\right\}$ is a normal subgroup of $G\rtimes H$ which is isomorphic to $G$, and $\overline{H}=\left\{\left(1,h\right)|h\in H\right\}$ is also a subgroup of $G\rtimes H$ which is isomorphic to $H$. Further, $\overline{G}\cap\overline{H}= 1$ and $\overline{GH}= G\rtimes H$. Identifying $\overline{G}$ with $G$ and $H$ with $\overline{H}$, one may view $G$ and $H$ as subgroups of $G\rtimes H$, and we will often write the element $\left(g,h\right)$ of $G\rtimes H$ as $gh$. 

	We are going to utilize the 'function notation' in describing wreath product with its particular product action. If we have that $G$ is a group, let $\Delta$ be a set, and let $H$ be the subgroup of $ Sym\Delta$. Since our focus will be on Cartesian decompositions, which, by definition, are finite, we shall throughout assume that the set $\Delta$ is finite. Take $B:=Func\left(\Delta,G\right)$, be a collection of all functions defined from $\Delta$ into $G$. Since $B$ being a given group defined based on the pointwise multiplication of elements of $B$. It has subgroups $G_{\delta}$, for $\delta\in\Delta$, defined 
$$
G=\left\{f\in Func\left(\Delta,G\right)\ |\ \delta'f=1\ for\ all\ \delta'\in\Delta\backslash\left\{\delta\right\}\right\}
$$
 each $G_{\delta}$ is isomorphic to $G$. Additionally, $B$ is known to be isomorphic with the direct product of these $|\Delta|$ duplicates of the group $G$, and the mapping $\sigma_{\delta}: f\mapsto f_{\delta}$ where $\delta'f_{\delta}=\left\{\begin{array}{l}
\delta f\ if\ \ \delta'=\delta\\
 1\ \ if\ \delta'\neq\delta
\end{array}\right.$ is defined as the natural projection mapping $B \rightarrow G_{\delta}$.

We then give a defininition of a group homomorphism, namely $\tau$, defined from the group $H$ to $Aut\left(B\right)$: Let $h\in H$ and $f\in B$. We take $f\left(h\tau\right)$ to be a function that is mapping  $f\left(h\tau\right):\delta\mapsto\delta h^{-1}f$. 

It is routine to check that $\tau$ is indeed a homomorphism. The wreath product $GwrH$ of $G$by $H$ is known in a general sense as semidirect product of $B$ by $H,$ i.e. $B\rtimes H$ with respect to the given homomorphism defined as $\tau$, and the subgroup $B$ of $GwrH$ is known as the {\it base group} of the $GwrH$, and group $H$ known as the {\it top group}. As the two components of a semidirect product are considered subgroups of the semidirect product, the base group $B$ and then the top group $H$ can also be considered subgroups of the wreath product and, in this way, $B$ becomes a normal subgroup of $GwrH$. Considering $H$ as a known subgroup of $GwrH$, the conjugation action of $H$ on $B$ will be induced by $\tau$ in $f\left(h\tau\right):\delta\mapsto\delta h^{-1}f$, and so we obtained:
$$
\left(\delta h^{-1}\right)f=\delta\left(f^{h}\right)\text{  }for\ all\ h\in H,\ f\in Func\left(\Delta,G\right),\ for\ all\ \delta\in\Delta.
$$
The wreath product $GwrH$ has a natural action by conjugation on the set of subgroups $G_{\delta}$ of its base group.

{\bf Remark 2.4}[9]{\bf :} Let us have a closer look at the special case $\Delta=\{1,2,\cdots,k\}$. Then, the wreath product $GwrH$ can be described using tuples, instead of functions. For $ i\in\Delta$, the image if can be written as $f_{i}$, and $f$ can be given by the $k$-tuple $\left(f_{1},f_{2},\cdots,f_{k}\right)$, If $h\in H$ and $\left(g_{1},g_{2},\cdots,g_{k}\right)\in B$ then the conjugate action of $h$ on $\left(g_{1},g_{2},\cdots,g_{k}\right)$ is known as
$$
\left(g_{1},g_{2},\cdots,g_{k}\right)^{h}=\left(g_{1h^{-1}},g_{2h^{-1}},\cdots,g_{kh^{-1}}\right)
$$
Hence $H$ permutes the coordinates of the elements of $B$.

{\bf Definition 2.5(Product action of Wreath Products} [1,9]{\bf )}

Given a product action of a wreath product $GwrH$ defined on the set of finctions $\Pi=Func(\Delta,\Gamma)$ to be: For $ G\leq Sym\Gamma$ and $ H\leq Sym\Delta$, let $f\in Func\left(\Delta,G\right)$ and with $h$ being an elemt of $H$ and fix

 $g=fh$. Let $\phi\in\Pi$, then we define the function $\phi g$ that maps the element $\delta\in\Delta$ to
$$
\delta\left(\phi g\right)=\left(\delta h^{-1}\phi\right)\left(\delta h^{-1}f\right)
$$
Observe the element $\delta h^{-1}\phi\in\Gamma$, with $\delta h^{-1}f\in Sym\Gamma$ giving $\left(\delta h^{-1}\phi\right)\left(\delta h^{-1}f\right)\in\Gamma$, and $\phi g\in Func\left(\Delta,\Gamma\right)=\Pi$, as expected.

Since $\Delta$ is finite, it is good to express the product action of it's given wreath product in the form of a coordinate notation. Suppose that $\Delta =\{1,2,\ldots,k\}$, and view $Func\left(\Delta,G\right)$ and $Func\left(\Delta,\Gamma\right)$ as $G^{k}$ and $\Gamma^{k}$, respectively, as in Remark 2.4. Then for $\left(\gamma_{1},\gamma_{2},\cdots,\gamma_{k}\right)\in\Gamma^{k},\ \left(g_{1},g_{2},\cdots,g_{k}\right)h\in GwrH$  we have
$$
\left(\gamma_{1},\gamma_{2},\cdots,\gamma_{k}\right)\left(\left(g_{1},g_{2},\cdots,g_{k}\right)h\right)=\left(\gamma_{1h^{-1}}g_{2h^{-1}},\gamma_{2h^{-1}}g_{2h^{-1}},\cdots,\gamma_{1h^{-1}}g_{kh^{-1}}\right)
$$
To have a deeper understanding of the subgroup(s) of wreath products of groups we entreat the concept of Cartesian decomposition. 

Looking at the set $\Pi = Func\left(\Delta,\Gamma\right)$, and $\forall \delta\in\Delta$, we give the defininition of a partition $\Gamma_{\delta}$ of the set $\Pi$ tobe: Take
\begin{center}
$\Gamma_{\delta} = \left\{\gamma_{\delta}|\ \forall\gamma\in\Gamma\right\},\ where \gamma_{\delta} :=\left\{\psi\in\Pi\ |\ \delta\psi=\gamma\right\}$.
\end{center}
We can easily check that $\Gamma_{\delta}$ is certainly a partition of of the set $\Pi$. The representation shows two significant facts. 

	First, the mapping $\delta\ \mapsto\ \Gamma_{\delta}$ is defined as a ono-to-one correspondence between $\Delta$ and $\{\Gamma_{\delta}|\delta\in\Delta\}$. 

	Second, let $\delta\in\Delta$, be a fixed element, the mapping $\gamma\mapsto\gamma_{\delta}$ defined is a given one-to-one correspondence between $\Gamma$ and $\Gamma_{\delta}$. Let $\gamma\in\Gamma$ and $\delta\in\Delta$, then the element $\gamma_{\delta}\in\Gamma_{\delta}$ is well thought-out to be 'a copy' of $\gamma$ in the set $\Gamma_{\delta}$, and it is known as the $\gamma$-part of $\Gamma_{\delta}$.

	The Cartesian product $\displaystyle \prod_{\delta\in\Delta}\Gamma_{\delta}$ is given to be a one-to-one correspondence with the original set $\Pi$: taking $\gamma_{\delta}\in\Gamma_{\delta}$, one for all element $\delta\in\Delta,$ the intersection $\displaystyle \bigcap_{\delta\in\Delta}\gamma_{\delta}$ comprises of single element of $\Pi$, namely the map that takes each $\delta$ to the element $\gamma\in\Gamma$ that corresponds to $\gamma_{\delta}$. This gives a one-to-one corespondence from the Cartesian product $\displaystyle \prod_{\delta\in\Delta}\Gamma_{\delta}$ to the set of functions $\Pi$. Then, the set
$$
\epsilon\ =\ \{\Gamma_{\delta}\ |\ \forall\ \delta\in\Delta\}
$$
is a Cartesian decomposition of the set $\Pi$. Precisely, is a set of partitions seen as the sets of  natural Cartesian decomposition of the set $\Pi$. Since $ Sym\Gamma wrSym\Delta$ is also a group that acts on the set $\Pi$, and the given action of $ Sym\Gamma wrSym\Delta$ is being strectched to subsets of the set $\Pi$, subsets of subsets, etc. To be specific, we look at the action of the group $ Sym\Gamma wrSym\Delta$ on the particular sets of partitions of the set $\Pi$. We will observe that $\{\Gamma_{\delta} | \delta\in \Delta\}$ is invariant under the action. The normal product action of the group $ Sym\Gamma wrSym\Delta$ on $\displaystyle \prod_{\delta\in\Delta}\Gamma_{\delta}$ is known to be permutationally isomorphic to the action that is defined on the set $\Pi$, and so the stabiliser in the permutation group $ Sym\Pi$ of the Cartesian decomposition is the permutation group $ Sym\Gamma wrSym\Delta$.

	If $X$ is a given subgroup of $ Sym\Gamma wrSym\Delta$ in its product action on the set of functions $\Pi=Func(\Delta,\Gamma)$, and $\Delta =\{1,\cdots,k\}$, then we identify $\Pi$ with the set $\Gamma^{k}$ of ordered $k$-tuples of the elements of the set $\Gamma$, and in this situation, subgroups of $ Sym\Gamma wrSym\Delta$ ascend as automorphism groups of different types of graph products (see (Praeger and Schneider, 2018b)), as groups of automorphism of some given codes of length say $k$ on the alphabet $\Gamma$, seen as subsets of $\Gamma^{k}.$

	Then following theorem and its proof can be found in Praeger and Schneider, 2018b, Theorem 5.13.

{\bf Theorem 2.6: (Wreath Embedding Theorem} [9]{\bf )} Suppose that $X$ is a given permutation group on a set $\Omega$ preserving a homogeneous Cartesian decomposition $\epsilon = \{\Gamma_{\delta}|\delta\in\Delta\}$  of $\Omega$, and let $\Gamma\in\epsilon$. Then there is a permutational isomorphism that maps $X$ to a subgroup of $ Sym\Gamma wrSym\Delta$ with its product action on $Func(\Delta,\Gamma)$, and maps $\epsilon$ to the natural Cartesian decomposition of $\epsilon$ defined above.

{\bf 3. MAIN RESULTS}

Now we are in a better position to prove our results. Embedding a wreath product in it's product action is equivalent to proving that a Cartesian decomposition is preserved. 

{\bf Definition 3.1}[3]{\bf :} The $G$-action is said to be transitive if $\Omega$ is a $G$-orbit; that is, for all $\alpha,\beta\in\Omega$ there is a given element $g\in G$ such that $\alpha g=\beta$. If $G$ is not transitive, then it is known as {\it intransitive}. A permutation group is known as {\it semiregular} if all its point stabilisers are trivial. A permutation group is {\it regular} if it is transitive and semiregular.

{\bf Proposition 3.2:} Let $T=S^{k}$, for some group $S$, let $S$ act regularly on $\Gamma$ such that $|\Gamma|=|S|$,  then $Sym\Gamma wrS_{k}$ in its natural action acts on $\Omega=\Gamma^{k}$, and the permutation representation of $T$ is embeddable in $Sym\Gamma wrS_{k}$ acting regularly on  $\Omega=\Gamma^{k}$.

{\bf Proof:} We suppose that $S$ act regularly on $\Gamma$ such that $|\Gamma|=|S|,\ T=S^{k}$, for some group $S$ acts regularly on $\Omega=\Gamma^{k}$. Let $\Omega=\Gamma^{k}$ and suppose there exist a given subgroup $W$ of $ Sym\Omega$ that is permutationally isomorphic to $Sym\Gamma wrS_{k},$ with $|\Gamma|\geq 2$ and $k\geq 2$, then $W=Sym\Gamma wrS_{k}$.

The normal subgroup $N=\left(Sym\Gamma\right)^{k}$ is known as the base group of $W$ and $H\cong S_{k}$ is known as the top group. The product action of $W$ on $\Omega=\Gamma^{k}$ is defined by
$$
\left(\gamma_{1},\cdots,\gamma_{k}\right)^{xh}=\left(\gamma_{1h^{-1}}^{x_{1h^{-1}}},\cdots,\gamma_{kh^{-1}}^{x_{kh^{-1}}}\right)
$$
for all $\left(\gamma_{1},\cdots,\gamma_{k}\right)\in\Omega, x=\left(x_{1},\cdots,x_{k}\right)\in N$ and $h\in H$, where the image of $\gamma\in\Gamma$ under $ y\in Sym\Gamma$ is $\gamma^{y}.\ W$ is obviously transitive on $\Omega=\Gamma^{k}$.

The Cartesian decomposition corresponding to the identity map on $\Omega$ is 
$$
\epsilon=\{\Gamma_{1},\cdots,\Gamma_{k}\}
$$
Where $\Gamma_{i}$ is the partition of $\Omega=\Gamma^{k}$ into disjoint subsets according to the $i^{th}$ coordinate of a point in $\Omega=\Gamma^{k}$, that is to say, the parts of $\Gamma_{i}$ are indexed by $\Gamma$ and the $\gamma$-part is the set of all points $(\gamma_{1},\cdots,\gamma_{k})$ with $\gamma_{i}=\gamma$. Thus $|\Gamma_{i}|=|\Gamma|$ for all $i$.

Thus $\epsilon$ is homogenous. Also each element $xh\in W$ maps the partition $\Gamma_{i}$ to the partition $\Gamma_{ih}$. Thus $W$ preserves the Cartesian decomposition $\epsilon$.

Also $W$ permutes the partitions $\Gamma_{i}$ transitively. Thus the permutation representation of $T$ is isomorphic to a subgroup of $Sym\Gamma wrS_{k}$.

{\bf Definition 3.3: (Diagonal group} D(T,m) $[3]${\bf )} Suppose that $G$ is a group with order $|G|>1$, and $n\geq 0$  positive integer. Let $\delta\left(G,n+1\right)$ be the diagonal subgroup $\left\{\left(g,g,\cdots,g\right)|g\in G\right\}$ of $G^{n+1}$. We select coset representatives for the element $\delta(G,n+1)$ in $G^{n+1}$. A suitable selection is to figure out the direct factors of $G^{n+1}$ as $G_{0};G_{1};\cdots;G_{n},$ and employ the representatives of the form $\left(1,g_{1},g_{2},\cdots,g_{n}\right)$ where $g_{i}\in G_{i}$. and let $\Omega$ denote the collection of all such symbols. Then, $\Omega$ is bijective by means of $G^{n}$.

We are now going to designate the action of $D(G,n)$ as:

(a) Let $1\leq i\leq n$, the factor $G_{i}$ acts by right multiplication on symbols in the $i^{th}$ position in the elements of the set $\Omega$.

(b) $G_{0}$ is acting by simultaneous left multiplication of all the coordinates by the inverse. Since, for $x\in G_{0},\ x$ maps the coset containing $\left(1,g_{1},g_{2},\cdots,g_{n}\right)$ to the coset containing $\left(x,g_{1},g_{2},\cdots,g_{n}\right)$, which is equal to the coset containing $\left(1,x^{-1}g_{1},x^{-1}g_{2},\cdots,x^{-1}g_{n}\right)$ Automorphisms of $G$ also acts simultaneously on each of the coordinates; nevertheless the inner automorphisms are recognized with the action of elements in the diagonal subgroup $\delta\left(G,n+1\right)$ (the element $(\left(x,x,x,\cdots,x\right)$ maps the coset containing $\left(1,g_{1},g_{2},\cdots,g_{n}\right)$  to the coset containing $\left(x,g_{1}x,g_{2}x,\cdots,g_{n}x\right)$, which is equal to the coset containing $\left(1,x^{-1}g_{1}x,x^{-1}g_{2}x,\cdots,x^{-1}g_{n}x\right).$

(c) Elements of the symmetric group $S_{n}$ (fixing coordinate $0$) also acts by permuting the coordinates in elements of $\Omega$.

(d) Look at the element of $S_{n+1}$ which transposes the coordinates $0$ and $1$. It is mapping the coset containing $\left(1,g_{1},g_{2},\cdots,g_{n}\right)$ to the coset containing $\left(g_{1},1,g_{2},\cdots,g_{n}\right)$ ,that also contains $\left(1,g_{1}^{-1},g_{1}^{-1}g_{2},\cdots,g_{1}^{-1}g_{n}\right)$. So the action of the given transposition is
$$
\left(1,g_{1},g_{2},\cdots,g_{n}\right)\mapsto\left(1,g_{1}^{-1},g_{1}^{-1}g_{2},\cdots,g_{1}^{-1}g_{n}\right).\ 
$$
({\it e}) Now $S_{n}$ and the transposition generates $S_{n+1}.$

{\bf Proposition 3.4:} Let $D\left(G,k\right)$ be a diagonal group where $G$ is a group such that $ G\leq Sym\Gamma$ and positive integer $k\geq 2$, then $D\left(G,k\right)$ is embeddable in the wreath product $Sym\Gamma wrS_{k}$ where the wreath product acts on $\Omega=\Gamma^{k}$ and $\Gamma\geq 2$ is a non-empty set. 

{\bf Proof:} Suppose that $D\left(G,k\right)$ is a diagonal group where $G$ is a group and there is a positive integer $k\geq 2$. 

Now, $Sym\Gamma wrS_{k}$ acts naturally on the set $\Omega=\Gamma^{k}$ in its product action and is defined by
$$
\left(\gamma_{1},\cdots,\gamma_{k}\right)^{xh}=\left(\gamma_{1h^{-1}}^{x_{1h^{-1}}},\cdots,\gamma_{kh^{-1}}^{x_{1h^{-1}}}\right)
$$
for all $\left(\gamma_{1},\cdots,\gamma_{k}\right)\in\Omega, x=\left(x_{1},\cdots,x_{k}\right)\in\left(Sym\Gamma\right)^{k}$ and $h\in S_{k}$ and the image of $\gamma\in\Gamma$ under $ y\in Sym\Gamma$ is $\gamma^{y}.\ Sym\Gamma wrS_{k}$ is obviously transitive on $\Omega=\Gamma^{k}$.

Now for each coset of the diagonal group $\delta(G,k)=\left\{\left(g,g,\cdots,g\right)|g\in G\right\}$ of $G^{k}$, there is a unique representative of the form $\left(1,g_{1},g_{2},\cdots,g_{k-1}\right)$ and define $\eta$ as
$$
\eta\left(\left(1,g_{1},g_{2},\cdots,g_{k-1}\right),\left(\left(g,g,...,g\right),\left(g,g,\cdots,g\right)\right)\right)=\left(1,g^{-1}g_{1}g,g^{-1}g_{2}g,\cdots,g^{-1}g_{k}g\right)
$$
which is the action on $\delta\left(G,k\right)$.

The Cartesian decomposition corresponding to the identity map on $\Omega$ is 
\begin{center}
$\epsilon=\{\Gamma_{1},\cdots,\Gamma_{k}\}$
\end{center}
Where $\Gamma_{i}$ is the partition of $\Omega=\Gamma^{k}$ into disjoint subsets according to the $i^{th}$ coordinate of a point in $\Omega=\Gamma^{k}$, that is to say, the parts of $\Gamma_{i}$ are indexed by $\Gamma$ and the $\gamma$-part is the set of all points $\left(\gamma_{1},\cdots,\gamma_{k}\right)$ with $\gamma_{1}=\gamma$. Thus $|\Gamma_{i}|=|\Gamma|$ for all $i$.

Thus $\epsilon$ is homogenous. Also each element of $Sym\Gamma wrS_{k}$ maps the partition $\Gamma_{i}$ to the partition $\Gamma_{ik}$. Thus $W$ preserves the Cartesian decomposition $\epsilon$. Thus $D\left(G,k\right)$ is isomorphic to a subgroup of $Sym\Gamma wrS_{k}.$

{\bf Definition 3.5:} Suppose that $G$ is a given group and $n$ a positive. Let $\Omega=G^{n}$ ; this will be the domain of a permutation, and its elements are written as $\left[x_{1},x_{2},\cdots,x_{n}\right]$, where $x_{1},x_{2},\cdots,x_{n}\in G$. The diagonal group $D\left(G,n\right)$ is generated by the following five types of permutations on $\Omega$:

(a) The group $G^{n}$ acting by right multiplication; so the element $\left(g_{1},g_{2},...,g_{n}\right)$ maps $\left[x_{1},x_{2},\cdots,x_{n}\right]$ to $\left[x_{1}g_{1},x_{2}g_{2},\cdots,x_{n}g_{n}\right]$. I will let $G_{i}$ be the $i^{th}$ factor of $G^{n}$, so that $G_{i}$ acts on the $i^{th}$ coordinate of elements of $\Omega$.

(b) The group $G$, acting by simultaneous left multiplication; so $g$ maps $\left[x_{1},x_{2},\cdots,x_{n}\right]$ to $\left[g^{-1}x_{1},g^{-1}x_{2},\cdots,g^{-1}x_{n}\right]$. It will denoted by $G_{0}$.

(c) The automorphism group of $G$, acting simultaneously on all coordinates.

(d) The symmetric group $S_{n}$, acting by permuting the coordinates.

(e) A permutation $\tau$ , defined by
\begin{center}
$\tau : \left[x_{1},x_{2},\cdots,x_{n}\right]\mapsto\left[x_{1}^{-1},x_{1}^{-1}x_{2},\cdots,x_{1}^{-1}x_{n}\right]$. 
\end{center}
{\bf Proposition 3.6:} Let $D\left(G,n\right)$ be a diagonal group, where $G$ is a finite group and $n$ is a positive integer, and assume that $\Gamma=G$ with $G$ acting regularly by right multiplication, and identify $G$ with the corresponding subgroup of $ Sym\Gamma$. Suppose also that in $Aut\left(G\right)$ there is a subgroup $Out\left(G\right)$, which complements the group of inner automorphisms. Then $D\left(G,n\right)$ is embeddable into the wreath product $Sym\Gamma wrS_{n}$, acting naturally in product action on the set $\Omega=\Gamma^{n}$.

{\bf Proof:} Let $D\left(G,n\right)$ be a diagonal group for a finite group $G$ and n a positive integer, and let $\Gamma=G$ with $G$ acting by right multiplication. Now, $W=Sym\Gamma wrS_{n}$ acts naturally on the set $\Omega=\Gamma^{n}$ in its product action and is defined by
$$
\left(\gamma_{1},\cdots,\gamma_{n}\right)^{nh}=\left(\gamma_{1h^{-1}}^{x_{1h^{-1}}},\cdots,\gamma_{nh^{-1}}^{x_{nh^{-1}}}\right)
$$
for all $\left(\gamma_{1},\cdots,\gamma_{n}\right)\in\Omega,\ x=\left(x_{1},\cdots,x_{n}\right)\in\left(Sym\Gamma\right)^{n}$ and $h\in S_{n}$ and the image of $\gamma\in\Gamma$ under $ y\in Sym\Gamma$ is $\gamma^{y}.\ W=Sym\Gamma wrS_{n}$ is obviously transitive on $\Omega=\Gamma^{n}.\ \Omega$ is bijective with $G^{n}$.

Now, 

	$\mathbf{1}\mathbf{.}$We identified $G$ with a subgroup of $ Sym\Gamma$. So we have $M:=G^{n}$ as a subgroup of the base group $\left(Sym\Gamma\right)^{n}$ of $W$.

	$\mathbf{2}\mathbf{.}$ From the definition of the action of $W$ it is clear that the top group $S_{n}$ normalizes $M$. 

	$\mathbf{3}\mathbf{.}$ Also, $N_{Sym\Gamma}\left(G\right)$ is the holomorph of $G$, this is known as a semidirect product of $G$ and a group $A=Aut\left(G\right)$ and $A$ acts on $\Gamma=G$ naturally as automorphisms. 

	$\mathbf{4}\mathbf{.}$ By assumption $A$ has a subgroup  $O\subseteq Out\left(G\right)$ which complements the inner automorphism group of $G$.

	$\mathbf{5}\mathbf{.}$ Hence the normalizer of $M$ in the base group of $W$ contains a semidirect product $MO^{n}$, and we define $D$ as the diagonal subgroup of $O^{n}$, namely $D=\left\{\left(x,x,\cdots,x\right)\ |\ x\in O\right\}.$

	$\mathbf{6}\mathbf{.}$ Finally consider the group generated by $M,\ D$ and $S_{n}$,  this is a copy of $D\left(G,n\right)$ in $W$.

The Cartesian decomposition corresponding to the identity map on $\Omega=\Gamma^{n}$ is 
$$
\epsilon=\{\Gamma_{1},\cdots,\Gamma_{n}\}
$$
Where $\Gamma_{i}$ is the partition of $\Omega=\Gamma^{n}$ into disjoint subsets according to the $i^{th}$ coordinate of a point in $\Omega=\Gamma^{n}$, that is to say, the parts of $\Gamma_{i}$ are indexed by $\Gamma$ and the $\gamma$-part is the set of all points $\left(\gamma_{1},\cdots,\gamma_{n}\right)$ with $\gamma_{i}=\gamma$. Thus $|\Gamma_{i}|=|\Gamma|$ for all $i$.

Thus $\epsilon$ is homogenous. Also each element of $Sym\Gamma wrS_{n}$ maps the partition $\Gamma_{i}$ to the partition $\Gamma_{ik}$. Thus $W$ preserves the Cartesian decomposition $\epsilon$. Now since the automorphism group of $G$ acts simultaneously on all coordinates, thus $W=Sym\Gamma wrS_{n}$ preserves the Cartesian decomposition $\epsilon$ and we conclude that $D\left(G,n\right)$ is embeddable into the wreath product $Sym\Gamma wrS_{n}.$

	{\bf 4. CONCLUSION}

We proved that a group that is acting regularly on a set and a diagonal group acting on a set in product action are embeddedable into wreath products in such actions using the Cartesian decomposition.

	{\bf 5. ACKNOWLEDGEMENT}

We would like to sincerely appreciate Prof. Peter J. Cameron of University of St Andrews, UK for his suggestion and also Prof. Cheryl E. Praeger of University of Western Australia, Australia for proof reading, insightful comments and suggestion. Also our appreciation goes to Prof. Csaba Schneider of Universidade Federal De Minas Gerais, Brazil for his insightful comment. 

	{\bf 6. REFRENCES}

[1] Praeger, C. E., and Schneider, C. Embedding Permustation Groups Into Wreath Products In Product Action. J. Aust. Math. Soc., (2012) 92, 127-136. DOI:10.1017/S1446788712000110

[2] Preager, C. E., and Schneider, C. Group Factorisations, Uniform Automorphisms, and Permutation Groups of Simple Diagonal Type. Israel Journal of Mathematics, (2018a) 228, 1001-1023. doi:DOI: 10.107/s11856-018-1790-1

[3] Dixon, J. D., and Mortimer, B. Permutation Groups. Berlin, New York: Springer. (1996).

[4] Sury, B. Wreath Products, Sylow's Theorem and Fermat's Little Theorem. European Journal of Pure and Applied Mathematics, (2010) 3(1), 13-15.

[5]  Schneider, C. Wreath products in permutation group theory. Retrieved September 3, 2019, from www.mat.ufmg.br/~csaba/pdf/marienheide.pdf (2017)

[6] Suleiman, E., and Audu, M. S. Some Embedment of Groups into Wreath Products. International Journal of Algebra and Statistics, (2020) 9(1-2), 1-10. doi:DOI: 10.20454/ijas.2020.1630

[7] Suleiman, E., and Audu, M. S. On Some Embedment of Groups into Wreath Products. Advances in Pure Mathematics, (2021) 11, 109-120. doi:https://doi.org/10.4236/apm.2021.112007

[8] Bailey, R. A., Cameron, P. J., Praeger, C. E., and Schneider, C. (2020). The geometry of diagonal groups. arXiv:2007.10726v1 [math.GR], (2020)1-61.

[9] Praeger, C. E., and Schneider, C. Permutation Groups and Cartesian Decompositions. (L. M. 449, Ed.) United Kingdom: Cambridge University Press. doi:DOI: 10.1017/9781139194006 (2018b)

[10] Praeger, C. E. The Inclusion Problem For Finite Primitive Permutation Groups. Proc. London Math. Soc., (1990) 60(3), 68-88.

\end{document}